\documentclass[a4paper,12pt]{amsart}\usepackage{amsfonts,amsmath,amssymb,bbm}

\def\C{\Bbb{C}}
\def\k{\mathbbm{k}}

\def\P{\Bbb{P}}\def\Q{\Bbb{Q}}\def\R{\Bbb{R}}\def\bX{{\bar{X}}}\def\Z{\Bbb{Z}}

\def\li{\ \\ $\bullet$ }\def\di{\partial}

\def\suml{\sum\limits}

\def\cupl{\mathop\cup\limits}\def\prodl{\mathop\prod\limits}

\newcommand{\quotients}[2]{{\footnotesize\left.\raisebox{0.4ex}{$#1$}\! /
\!\raisebox{-0.4ex}{$#2$}\right.}}

\def\tX{{\tilde{X}}}

\def\hx{{\hat{x}}}

\def\al{\alpha}\def\Ga{\Gamma}\def\de{\delta}\def\De{\Delta}
\def\ep{\epsilon}
\def\la{\lambda}\def\si{\sigma}

\def\cK{{\mathcal K}}\def\cL{{\mathcal L}}\def\cN{{\mathcal
N}}\def\cO{\mathcal O}

\def\cY{{\mathcal Y}}

\def\norm{\stackrel{\nu}{\to}}

\def\u0{\underline{0}}
\def\uf{\underline{f}}
\def\uk{{\underline{k}}}

\def\ux{{\underline{x}}}

\def\empty{\varnothing}

\newcommand{\ber}{\begin{array}{l}}\newcommand{\eer}{\end{array}}
\newcommand{\bpm}{\begin{pmatrix}}\newcommand{\epm}{\end{pmatrix}}
\newcommand{\bM}{\begin{matrix}}\newcommand{\eM}{\end{matrix}}
\newcommand{\bee}{\begin{enumerate}}\newcommand{\eee}{\end{enumerate}}

\def\wrt{with respect to
}\def\sset{\subset}\def\sseteq{\subseteq}\def\smin{\setminus}

\def\cf{C_{n,r}}\newcommand{\stir}[2]{\genfrac{\{}{\}}{0pt}{0}{#1}{#2}}

\def\Nnd{Newton-non-degenerate}\def\ND{Newton diagram}
\def\wrt{with respect to }

\def\bull{\vrule height .9ex width .9ex depth -.1ex }

\newcommand{\beq}{\begin{equation}}\newcommand{\eeq}{\end{equation}}
\newcommand{\beqm}{\begin{multline}}\newcommand{\eeqm}{\end{multline}}
\newtheorem{Lemma}{Lemma}[section]\newcommand{\bel}{\begin{Lemma}}\newcommand{\eel}{\end{Lemma}}
\newtheorem{Example}[Lemma]{Example}\newcommand{\bex}{\begin{Example}\rm}\newcommand{\eex}{\end{Example}}
\newtheorem{Proposition}[Lemma]{Proposition}\newcommand{\bprop}{\begin{Proposition}}\newcommand{\eprop}{\end{Proposition}}
\newtheorem{Definition-Proposition}[Lemma]{Definition-Proposition}

\def\bpr{~\\{\em Proof.\ }}\newcommand{\epr}{$\bull$\\}
\newtheorem{Theorem}[Lemma]{Theorem}\newcommand{\bthe}{\begin{Theorem}}\newcommand{\ethe}{\end{Theorem}}
\newtheorem{Definition}[Lemma]{Definition}\newcommand{\bed}{\begin{Definition}}\newcommand{\eed}{\end{Definition}}
\newtheorem{Remark}[Lemma]{Remark}\newcommand{\beR}{\begin{Remark}\rm}\newcommand{\eeR}{\end{Remark}}
\newtheorem{Corollary}[Lemma]{Corollary}\newcommand{\bcor}{\begin{Corollary}}\newcommand{\ecor}{\end{Corollary}}
\newcommand{\bet}{\begin{tabular}{cccccccc}}\newcommand{\eet}{\end{tabular}}

\DeclareMathOperator{\const}{const}
\DeclareMathOperator{\mult}{mult}
\DeclareMathOperator{\Vol}{Vol}
\DeclareMathOperator{\Covol}{Covol}
\DeclareMathOperator{\Cone}{Cone}
\DeclareMathOperator{\Conv}{Conv}
\DeclareMathOperator{\Supp}{Supp}
\DeclareMathOperator{\Span}{Span}
\DeclareMathOperator{\rank}{rank}
\DeclareMathOperator{\lt}{l.t.}
\DeclareMathOperator{\lot}{l.o.t.}

\newcommand{\bin}[2]{\binom{#1}{#2}}

\title[]{D\MakeLowercase{urfee-type bound for some non-degenerate complete intersection
singularities}}
\author[]{D\MakeLowercase{mitry} K\MakeLowercase{erner} \MakeLowercase{and}
A\MakeLowercase{ndr\'as} N\MakeLowercase{\'emethi}}
\address{\noindent Department of Mathematics, Ben Gurion University of the Negev, P.O.B. 653, Be'er Sheva 84105,
Israel, dmitry.kerner@gmail.com,
\newline
R\'enyi Institute of Mathematics, Budapest, Re\'altanoda u. 13--15, 1053, Hungary,
nemethi@renyi.hu}
\date{\today}
\thanks{D.K. was supported by the grant FP7-People-MCA-CIG, 334347.
\\ \ A.N. was supported by  OTKA Grant 100796}

\begin{document}
\begin{abstract}
The Milnor number, $\mu(X,0)$, and the singularity genus, $p_g(X,0)$, are  fundamental invariants of isolated hypersurface singularities (more generally, of local complete intersections). The long standing Durfee conjecture (and its generalization) predicted the inequality $\mu(X,0)\ge (n+1)!p_g(X,0)$, here $n=\dim(X,0)$. Recently we have constructed counterexamples, proposed a corrected bound and verified it for the homogeneous complete intersections.

In the current paper we treat the case of germs with \Nnd\ principal part when the \ND s are ``large enough", i.e. they are large multiples of some other diagrams. In the case of local complete intersections we prove the corrected inequality, while in the hypersurface case we prove an even stronger inequality.
\end{abstract}
\maketitle \setcounter{secnumdepth}{6} \setcounter{tocdepth}{1}

\section{Introduction}
\subsection{} Let $(X,0)\sset(\C^{n+r},0)$ be the germ of an isolated analytic complete intersection
singularity of dimension $n$.
 The Milnor number and the singularity genus are fundamental local invariants. They can be
 defined as the defects of the corresponding global invariants.
Indeed, by the finite determinacy, we can assume $(X,0)$ to be an algebraic germ: let $X$ be a
 representative of $(X,0)$. Take some projective compactification $X\sset\bX$, assume no
 other singularities are added,
  i.e. $\bX\smin X$ is smooth.
 Take (one of) its resolution, $\tX\to\bX$ and (one of) its smoothing $\bX_\ep$.
 Then
\beq\label{Eq.Def.of.Mu.Pg}
\mu(X,0):=(-1)^n\Big(\chi_{\text{top}}(\bX_\ep)-\chi_{\text{top}}(\bX)\Big),\quad
p_g(X,0):=(-1)^n\Big(\chi_{\text{an}}(\cO_\bX)-\chi_{\text{an}}(\cO_\tX)\Big).
\eeq
(Here $\chi_{\text{top}}$ is the topological Euler characteristic, while $\chi_{\text{an}}(\cO)$ is the
analytic Euler characteristic of the structure sheaf.)
These invariants do not depend on the choice of the resolution/smoothing/compactification, they
are totally determined by the local analytic geometry of the germ $(X,0)\sset(\C^{n+r},0)$.
(In fact, in the hypersurface case, $r=1$, $p_g$ is even preserved in the $\mu=\const$ deformations,
\cite[pg.115]{AGLV}.)

The relation between the Milnor number and the singularity genus has been investigated for long
time. For example, in the case of curves $p_g$ coincides with the classical $\de$-invariant of
the singularity.
 Then one has the relation $\de=\frac{\mu+r-1}{2}$, \cite{Buchweitz-Greuel-1980}, where $r$ is
 the number of local branches at the singular point.

In \cite{Durfee1978} the inequality $\mu\ge 6p_g$ was conjectured for surface singularities that are
 isolated complete intersections. In \cite{Kerner-Nemethi.Durfee1},
 \cite{Kerner-Nemethi.Durfee2}
 we disproved this initial inequality and proposed a modified inequality for isolated complete
 intersections $(X,0)\sset(\C^{n+r},0)$
 (of dimension $n>2$ and codimension $r$): $\mu\ge\cf p_g$. (For $n=2$ the only possible universal bound is $\mu\ge
 4p_g$.)
 We proved the new inequality for homogeneous complete intersections.

 The combinatorial coefficient $\cf$ is defined by
 $\cf:=\frac{\bin{n+r-1}{n}(n+r)!}{\stir{n+r}{r}r!}$.
 Here $\stir{n+r}{r}$ is the Stirling number of the second kind.
 For more details see \S\ref{Sec.Background.Coefficient.Cnr}, now we only quote the basic
 property:
\beq
(n+1)!=C_{n,1}> C_{n,2}> \cdots > C_{n,r} >\cdots > \lim_{r\to\infty}\cf=2^n.
\eeq

For the history and the list of other (partial) verifications see
\cite{Kerner-Nemethi.Durfee2}. For the relevant notions from Singularity Theory see \cite{AGLV}, \cite{Dimca92},
\cite{Looijenga-book}, \cite{Oka-book}.

This paper is the continuation of \cite{Kerner-Nemethi.Durfee1} and \cite{Kerner-Nemethi.Durfee2}.
We verify the corrected bound for several additional classes of
singularities.

\subsection{}
Our first main result is the following:
\bthe\label{Thm.Bound.for.Nnd.C.I.}
Consider an isolated complete intersection singularity,
$(X,0)=\{f_1=\cdots=f_r=0\}\sset(\C^{n+r},0)$.
Suppose either $r=1$, $n\ge2$, or $r>1$, $n>2$. Suppose the tuple $(f_1,\dots,f_r)$ is \Nnd\
\wrt\ the diagrams $(\Ga_1,\dots,\Ga_r)$. Suppose
  that all the diagrams are convenient and `large enough'. (Namely, for $i=1,\dots,r$:
  $\Ga_i=d_i\tilde{\Ga}_i$, where $1\ll d_i\in\Q$
  and $\tilde{\Ga}_i$ are some other fixed Newton diagrams.)
Then $\mu(X,0)> \cf p_g(X,0)$.

Further, the bound is asymptotically sharp (i.e. $\frac{\mu}{p_g}\to\cf$ for
$\max\{d_i\}_i\to\infty$) iff $\Ga_1=\cdots=\Ga_r$.
\ethe

\beR
1.  For \Nnd\ singularities the Milnor number  and the singularity genus are determined
combinatorially by the \ND s.
 Therefore in this case the proof of the inequality consists of a lattice point count and its
 comparison to the volume(s)
  of the bodies under the \ND s.

2. Even with this reduction to combinatorics, the proof is not straightforward. It is heavily
based on an `inequality of averages',
a  Fortuin--Kasteleyn--Ginibre-type result, which we prove separately in
\cite{Kerner-Nemethi.Combinatorial.Lemma}.

3. In \cite{Kerner-Nemethi.Durfee2} we have considered isolated complete intersections, when all
$\{f_i\}$ are homogeneous.
  In that case we proved that the bound is asymptotically sharp
 precisely when all the multiplicities coincide. Therefore, our present sharpness statement
 (`the bound is sharp  iff $\Ga_1=\cdots=\Ga_r$')
  is the natural extension of this fact.

4. Recall that $p_g$ is defined for singularities over any algebraically closed field $\k$ of
 zero characteristic.
 The Milnor number is a topological invariant, but in some cases it can be defined also
 for singularities over $\k$, and it satisfies the usual properties of the `classical'
  Milnor number.
Our proof is purely combinatorial, it does not use any complex topology. Therefore, if one
 defines a Milnor number over $\k$, with the usual properties (in particular if the results of
 Kouchnirenko-Khovanskii-Bivia-Ausina
   hold),   then our proof holds over $\k$ as well.
\eeR

\subsection{} For hypersurface singularities which are \Nnd\ and have  large enough \ND, we prove in
\S\ref{Sec.Stronger.Bound.Hypersurfaces}
  a stronger inequality:
\bthe\label{Thm.Asympt.Bound.Hypersurfaces} Assume  $n>2$ and let
$(X,0)\sset(\C^{n+1},0)$ be the germ
of an isolated hypersurface singularity, non-degenerate with respect
to its \ND. Let $p=\mult(X,0)$ be its multiplicity. Suppose the \ND\ of $(X,0)$
 is `large enough', i.e. $\Ga_{(X,0)}=d\tilde{\Ga}$, where $0\ll d\in\Q$, while $\tilde{\Ga}$ is
 some other \ND.
\\1. Then
 $\mu(X,0)-\Big((p-1)^{n+1}-\frac{p!}{(p-n-1)!}\Big)\ge(n+1)!p_g(X,0)$.
\\2. If, moreover, the projectivized tangent cone, $\P T_{(X,0)}$, has at most isolated singularities,
 with total Milnor number $\mu(\P T_{(X,0)})$, then
\[\mu(X,0)-\mu(\P T_{(X,0)})-\Big((p-1)^{n+1}-\frac{p!}{(p-n-1)!}\Big)\ge(n+1)!p_g(X,0).\]
 Here the equality holds iff
$\P T_{(X,0)}$ is smooth, i.e. $(X,0)$ is a homogeneous isolated hypersurface singularity.
\ethe
If the projective hypersurface $\P T_{(X,0)}\sset\P^n$ has only isolated singularities then the total Milnor number, $\mu(\P T_{(X,0)})$, is the sum of the local Milnor numbers, in particular it is positive.
When $\P T_{(X,0)}$ has non-isolated singularities the total Milnor number $\mu(\P T_{(X,0)})$, defined in equation (\ref{Eq.Def.of.Mu.Pg}), can be negative.
Whenever this term is negative we will not consider it in the inequality.

\

Probably one can extend this type of stronger inequality to the complete intersections and prove:
\beq
\mu-\cf
p_g>\suml_{\uk\in\cK_{n,r}}\Vol_{n+r}\big((\Ga^+_1)^{k_1},\dots,(\Ga^+_r)^{k_r}\big)\Big((n+r)!-\cf\bin{n+r}{k_1,k_2,\dots,k_r}\Big).
\eeq
 Here the right hand side has often the order
of $\mu$. It vanishes when all the diagrams are
  proportional (in particular it vanishes for $r=1$). So, this right hand side cannot be seen when all $f_i$
  are ordinary multiple points (i.e. isolated homogeneous singularities) or in the case of hypersurface singularities.

\subsection{} As one sees above, for \Nnd\ singularities we always assume that the
diagram(s) is/are
`large enough'. As of now we could not prove the Durfee bound
 for an arbitrary \ND, even for \Nnd\ surface singularities in $(\C^3,0)$.
 In this case the combinatorial formulas are:
\beq
\mu=\mu(\Ga)=3!\Vol_3(\Ga)-2!\Vol_2(\Ga)+\Vol_1(\Ga)-1,\quad p_g=p_g(\Ga)=|\Ga^-\cap\Z^3_{>0}|,
\eeq
while the conjectural bound is $\mu(\Ga)\ge 6p_g(\Ga)$. (Recall that $\Vol_i$ denotes the
normalized $i$-dimensional lattice volume, as e.g. in \cite{Kouchnirenko76}.)

It is natural to try to extend this (purely combinatorial) bound to some more general class of lattice
polytopes. The situation is highly delicate as the following example shows.
\bex
Suppose instead of \ND s in $\R^3_{\ge0}$ one considers a generalized version: \ND s inside the
cone $\Cone_{x^{-1}_1x^{-1}_2x^{-1}_3}(x^m_1,x^m_2,x^m_3)$.
 (This cone is generated by the rays starting from $x^{-1}_1x^{-1}_2x^{-1}_3$ and passing
 through any of $x^m_1$, $x^m_2$, $x^m_3$.)
This means that we consider \Nnd\ hypersurfaces inside a toric variety with an isolated
singularity. Consider the \ND\ of the homogeneous $\Ga=\Conv(x^m_1,x^m_2,x^m_3)$.
 Then $\Ga^-=\Conv(x^{-1}_1x^{-1}_2x^{-1}_3,x^m_1,x^m_2,x^m_3)$ and its parameters are (see
 \S\ref{Sec.Background.Newton.Diagrams}):
\begin{multline}
\Vol_3(\Ga^-)=\frac{m^3}{6}+3\frac{m^2}{6},\quad \Vol_2(\Ga^-)=3\cdot\frac{m}{2},\quad
\Vol_1(\Ga^-)=3\cdot 1,
\quad
\\|\overset{\circ}{\Ga^-}\cap\Z^3|=\bin{m+2}{3},\quad
|\overset{\circ}{\Ga}\cap\Z^3|=\bin{m-1}{2}.
\end{multline}
Then the singularity invariants are (see \S\ref{Sec.Background.mu.pg.for.Nnd.C.I.}):
\beq
\mu(\Ga^-)=m^3+3m^2-3m+2,\quad\quad p_g(\Ga^-)=\bin{m+2}{3}+\bin{m-1}{2}.
\eeq
Therefore: $\mu(\Ga^-)-6p_g(\Ga^-)=-3m^2+4m-4<0$, i.e. for this $\Ga^-$ the inequality is
violated.
\eex
Therefore, when trying to prove the inequality in the ordinary case, $\Ga^-\sset\R^3$, we cannot
subdivide the body $\Ga^-$ into some suitable pieces and combine the total $\mu>6p_g$ from its
building blocks. Geometrically, this inequality cannot be proven by any local consideration of
the resolution of $(X,0)$, rather it depends on its global properties.

\section{Preliminaries}
\subsection{Some relevant combinatorics}
\subsubsection{Stirling numbers}\label{Sec.Stirling.Coefficients.2'nd.kind}
For any $n\geq 0$ and $r\geq 1$,
 the Stirling number of the second kind, $\stir{n+r}{r}$, is the number of (unordered)
 partitions of $n+r$ elements into $r$
 non-empty sets, see \cite[\S24.1.4, pg. 824]{Abramowitz-Stegun}. We record some of its basic
 properties:
\li $\stir{n}{1}=1$, $\stir{r}{r}=1$;
\li for $r>1$: $\stir{n+r}{r}\ge\stir{n+r-1}{r-1}$ and the equality occurs only for $n=0$;
\li the generating function for these numbers is
$(e^x-1)^r=r!\suml^\infty_{n=0}\stir{n+r}{r}\frac{x^{n+r}}{(n+r)!}$;
\li the explicit expansion:
$\stir{n+r}{r}=\frac{1}{r!}\suml^r_{j=0}(-1)^j\bin{r}{j}(r-j)^{n+r}$;
\li the asymptotic growth: $\stir{n+r}{r}\sim\frac{r^{2n}}{2^nn!}$, as $r\to\infty$;
\li the recurrence relation: $\stir{n+1+r}{r}=r\stir{n+r}{r}+\stir{n+r-1}{r-1}$;
\li for $n+r\ge r\ge j$ there is another recurrence relation:

\beq\label{Eq.Recurrence.relation.Stirling.coeff}
\bin{r}{j}\stir{n+r}{r}=\suml^{n+r-j}_{i=r-j}\bin{n+r}{i}\stir{n-i+j}{j}\stir{i+r-j}{r-j}.
\eeq

\subsubsection{The set of compositions}\label{Sec.Knr.first.time}
Denote by $\cK_{n,r}$ the set of the (ordered) compositions,
\beq
\cK_{n,r}:=\{\uk=(k_1,\ldots, k_r)\, :\,  k_i\geq 0 \ \mbox{for all $i$, and } \ \sum
_ik_i=n\}.
\eeq
This $\cK_{n,r}$ can be thought of as the lattice points of the simplex.
 Its  cardinality is $|\cK_{n,r}|=\binom{n+r-1}{n}$.

The permutation group on $r$ elements, $\Xi_r$, acts on $\cK_{n,r}$.  The quotient
$\quotients{\cK_{n,r}}{\Xi_r}$ is the set of partitions. (Recall that a partition is an {\em
unordered} composition.)
For convenience we put $\cK_{n,r}=\varnothing$ when $r\le0$ or $n<0$.

Suppose a set of objects is indexed by this set of compositions, $\{A_\uk\}_{\uk\in\cK_{n,r}}$.
 We often use the standard set-theoretic inclusion-exclusion formula:
\beq\label{Eq.Exclusion.Inclusion.Formula}
\suml_{\uk\in\cK_{n,r}}A_\uk-\suml^r_{i=1}\suml_{\substack{\uk\in\cK_{n,r}\\k_i=0}}A_\uk+
\suml_{1\le i_1<i_2\le r}\suml_{\substack{\uk\in\cK_{n,r}\\k_{i_1}=0=k_{i_2}}}A_\uk-\cdots=
\suml_{\substack{\uk\in\cK_{n,r}\\k_1,\dots,k_r>0}}A_\uk.
\eeq

\subsubsection{The coefficient $\cf$}\label{Sec.Background.Coefficient.Cnr}
Using these notions the coefficient $\cf$ is defined by
\beq\label{eq:1}
\cf:=\frac{\bin{n+r-1}{n}(n+r)!}{\stir{n+r}{r}r!}=
\frac{|\cK_{n,r}|}{\suml_{\uk\in\cK_{n,r}}\prodl^r_{i=1}\frac{1}{(k_i+1)!}}.
\eeq
The second equality of (\ref{eq:1}) follows from \cite[pages 176-178]{Jordan1965}.

We record some properties of $\cf$.
\\1. $C_{n,1}=(n+1)!$, $C_{n,2}=\frac{(n+2)!(n+1)}{2^{n+2}-2}$,
$C_{n,3}=\frac{\bin{n+2}{2}(n+3)!}{3^{n+3}-3\cdot 2^{n+3}+3}$,
 by direct computation.
\\2. $\lim_{r\to\infty}\cf=2^n$. The limit can be  computed using the asymptotical growth of
Stirling
 numbers, \S\ref{Sec.Stirling.Coefficients.2'nd.kind}.
 This gives: $\cf\sim 2^n\frac{(n+r-1)!(n+r)!}{(r-1)!r!r^{2n}}$
with limit $2^n$ as $r\to\infty$.
\\3. $C_{n,1}> C_{n,2}> \cdots > C_{n,r} >\cdots > \lim_{r\to\infty}\cf.$
 This is proved e.g. in \cite[Corollary 4.2]{Kerner-Nemethi.Durfee2}.
\\4. $\suml_{\uk\in\cK_{n,r}}\Big[(n+r)!-C_{n,r}\bin{n+r}{k_1+1,\dots,k_r+1}\Big]=0$.
 This follows immediately from equation (\ref{eq:1}).

\subsection{\ND s}\label{Sec.Background.Newton.Diagrams}
Let $f(x_1,\dots,x_N)=\suml_I a_I \ux^I$ be a power series, with complex coefficients. Consider
the support
 of its monomials, $\Supp(f):=\{I\in\Z^N_{\ge0}|\ a_I\neq0\}$. The Newton polyhedron is defined
 as
 the convex hull $\Ga^+=\Ga^+_f:=\Conv(\Supp(f)+\R^N_{\ge0})$.
The Newton polyhedron has compact faces and unbounded faces.
 The Newton diagram, $\Ga$, is the union of all the compact faces of $\Ga^+$.
We always assume that the diagram is `convenient', that is,
 $\Ga$ intersects all
 the coordinate axes (i.e. $f$ contains all the monomials $x^{m_1}_1,\dots,x^{m_N}_N$).

We  use the notation  $\Ga^-:=(\R^N_{\ge0}\smin\Ga^+)\cup\Ga$ for the
{\em part not above the diagram}.
Denote the set of lattice points on $\Ga$ by $\Ga\cap\Z^N$, similarly $\Ga^-\cap\Z^N$.
 The notation for the lattice points strictly below the diagram is $(\Ga^-\smin\Ga)\cap\Z^N$.

Let $\Vol_N(\Ga^-)$ be the (lattice) volume of $\Ga^-$ (with the normalization that the volume of
the unit cube is 1.)
 More generally, for any subset $I\sseteq\{1,\dots,N\}$ consider the corresponding coordinate
 plane $L_I=\Span(\{\hx_i\}_{i\in I})$.
 (Here $\hx_i$ is the unit vector along the $i$'th coordinate axis.)
 Define $\Ga^I=\Ga\cap L_I$ and $(\Ga^-)^I=\Ga^-\cap L_I$.
 Accordingly, for a fixed $j$, one has the sum of volumes of intersections with all
$j$-dimensional coordinate planes:
\beq
\Vol_j(\Ga^-):=\suml_{\substack{I\sset\{1,\dots,N\}\\|I|=j}}\Vol_j\Big((\Ga^-)^I\Big).
\eeq
In
particular, $\Vol_0(\Ga^-)=1$, while $\Vol_1(\Ga^-)$ is the total
lattice length of all the segments of the type
$\Conv(\vec{0},d_i\hx_i)$, where $\Ga\cap \Span(\hx_i)=d_i\hx_i$.

The diagram $\Ga$ consists of many faces. Each face has its
(lattice) volume inside the lattice it spans. Let $\Vol_{N-1}\Ga$ be
the total volume of $\Ga$, i.e. the sum of the volumes of the top
dimensional faces.

\subsection{Mixed covolumes and their convexity}\label{Sec.Background.Mixed.Covolumes}
Given a convenient Newton polyhedron, $\Ga^+\sset\R^N_{\ge0}$, consider its covolume,
$\Covol(\Ga^+):=\Vol_N(\R^N_{\ge0}\smin\Ga^+)$.
Given a collection of Newton polyhedra, $\{\Ga^+_i\}_i$, consider their scaled Minkowski sum,
$\la_1\Ga^+_1+\cdots+\la_r\Ga^+_r$.
The covolume of this sum is a polynomial in $\{\la_i\}$,
\cite[\S10]{Kaveh-Khovanskii-Conv.Bod}:
\beq
\Covol(\la_1\Ga^+_1+\cdots+\la_r\Ga^+_r)=
\suml_{\uk\in\cK_{N,r}}\bin{N}{k_1,\dots,k_r}\Covol\Big((\Ga^+_1)^{k_1},
\dots,(\Ga^+_r)^{k_r}\Big)(\prod^r_{i=1}\la^{k_i}_i).
\eeq
The mixed covolumes are the (positive) coefficients
$\Covol\Big((\Ga^+_1)^{k_1},\dots,(\Ga^+_r)^{k_r}\Big)$.
\\Here $\Covol\Big((\Ga^+_1)^{k_1},\dots,(\Ga^+_r)^{k_r}\Big)$ is a shorthand for
$\Covol\Big(\underbrace{\Ga^+_1,\dots,\Ga^+_1}_{k_1},\dots,\underbrace{\Ga^+_r,\dots,\Ga^+_r}_{k_r}\Big)$,
for $k_1+\cdots+k_r=N$.

We use the following basic properties of the mixed covolumes:
\li They are symmetric and multilinear:
\[
\Covol(\Ga^+_{11}+\Ga^+_{12},\Ga^+_2,\dots,\Ga^+_N)=\Covol(\Ga^+_{11},\Ga^+_2,\dots,\Ga^+_N)+
\Covol(\Ga^+_{12},\Ga^+_2,\dots,\Ga^+_N).
\]
\li For the diagrams  $\{\Ga_i=d_i\Ga\}_{i=1,\dots,r}$ one has
\begin{multline}
\Covol_N(\suml_i\la_i\Ga^+_i)=\Covol_N(\suml_i d_i\la_i\Ga^+)=(\suml_i d_i\la_i)^N
\Covol_N(\Ga^+)=\\=
\suml_{\uk\in\cK_{n,r}}\bin{N}{k_1,\dots, k_r}(\prodl^r_{i=1} (\la_i d_i)^{k_i})
\Covol_N(\Ga^+).
\end{multline}
\li Convexity property: $\Covol(\Ga^+_1,\Ga^+_2,\cdots,\Ga^+_N)^2\le
\Covol(\Ga^+_1,\Ga^+_1,\Ga^+_3,\cdots,\Ga^+_N)\Covol(\Ga^+_2,\Ga^+_2,\Ga^+_3,\cdots,\Ga^+_N)$,
\cite{Teissier1978},
\cite[Appendix]{Teissier2004}, \cite[Theorem 10.5]{Kaveh-Khovanskii-Mix.Mult}.
\li In the proof of theorem \ref{Thm.Bound.for.Nnd.C.I.} we use the following generalization of
the convexity of the mixed co-volumes:
\beqm\label{Eq.Ineq.of.mixed.Covolumes}
\Big(\suml_{\uk\in\cK_{n,r}}\bin{n+r}{k_1+1,\dots,k_r+1}\Big)\cdot
\suml_{\substack{\uk\in\cK_{n+r,r}\\k_1,\dots,k_r\ge1}}
\Covol_{n+r}\Big((\Ga^+_1)^{k_1},\dots,(\Ga^+_r)^{k_r}\Big)
\ge\\
\ge
\bin{n+r-1}{n}\cdot \suml_{\substack{\uk\in\cK_{n+r,r}\\k_1,\dots,k_r\ge1}}
\bin{n+r}{k_1+1,\dots,k_r+1}\Covol_{n+r}\Big((\Ga^+_1)^{k_1},\dots,(\Ga^+_r)^{k_r}\Big).
\end{multline}
This inequality is proved separately in \cite[\S4]{Kerner-Nemethi.Combinatorial.Lemma}.

\subsection{Non-degeneracy \wrt\ \ND s}\label{Sec.Background.Non-degeneracy}
The non-degeneracy notion was studied first for functions in \cite{Kouchnirenko76},
 then for complete intersections in  \cite{Khovanskii1978}). The material of this section is
 taken from \cite[\S3]{Bivia-Ausina2007}, see also \cite{Bivia-Ausina2004}.

Consider several power series, $g_1,\dots,g_r\in\C[[x_1,\dots,x_N]]$, for $r\le N$. Take
Minkowski sum of
their Newton polyhedra, $\Ga^+:=\Ga^+_1+\cdots+\Ga^+_r$.
Let $\si$ be a compact face of $\Ga^+$.
By \cite[Lemma 2.7]{Damon1989} and \cite[Lemma 3.4]{Bivia-Ausina2007} there exists the unique
set
 of compact faces, $\si_1\sset\Ga_1$, \dots, $\si_r\sset\Ga_r$ satisfying:
 $\si=\si_1+\cdots+\si_r$.

The part of $g_i$ supported on $\sigma_i$ will be  denoted by $g_i|_{\sigma_i}$.
\bed
1. The sequence $g_1,\dots,g_r$ satisfies the
 $(B_\si)$ condition if $\{g_1|_{\si_1}(x)=\cdots=g_r|_{\si_r}(x)=0\}\cap(\C^*)^N=\empty$.
\\2. The sequence $g_1,\dots,g_r$ is non-degenerate if it is a regular sequence
 (i.e. defines a subspace of codimension $r$) and satisfies the $(B_\si)$ condition for all the
 compact faces $\si$ of $\Ga^+$ of dimension $\dim(\si)\le r-1$.
\eed
To define the non-degeneracy of the map $f=\{f_1,\dots,f_r\}$ we need the notion of
non-degeneracy of modules.

For any ideal $J$ the Newton polyhedron is defined by
 $\Ga^+(J)=\Conv\Big(\cupl_{f\in J}(\Supp(f)+\R^N_{\ge0})\Big)$. The
 Newton diagram (or the diagram of exponents) is defined as in
 \S\ref{Sec.Background.Newton.Diagrams}.

Consider a submodule of a free module, $M\sset \C[[x_1,\dots,x_N]]^{\oplus r}$. Denote by $A_M$
its generating
 matrix, i.e. a $r\times s$ matrix with entries in $\C[[x_1,\dots,x_N]]$, whose columns generate
 the module.
 Denote by $M_i$ the ideal in $\C[[x_1,\dots,x_N]]$ generated by the entries of $i$'th row of
 $A_M$.
  (It does not depend on the choice of generators of the module.) The {\em Newton polyhedron of
  $M$} is
   defined to be $\Ga^+(M):=\Ga^+(M_1)+\cdots+\Ga^+(M_r)$.
   (Here each $M_i$ is an ideal and we use the definition of $\Ga^+(J)$ as above.
   In the case of one-row-matrix $M$ itself is an ideal.)
    For any compact face $\si$ of $\Ga^+(M)$ take
   its (unique) presentation $\si=\si_1+\cdots+\si_r$, $\si_i\sset\Ga^+(M_i)$, as above.
   Denote by $M|_\si$ the matrix of restrictions, its $i$'th row consists of the restrictions
   onto $\si_i$.
   (Note that all the restrictions are polynomials, not just power series.)
\bed
The module/matrix $M$ is called \Nnd\ if for any compact face $\si\sset\Ga^+(M)$ the following property holds:
\[\{x\in\C^N:\ \rank(M|_\si(x))\le r\}\cap(\C^*)^N=\empty.\]
\eed
Finally, for a map $\uf=(f_1,\dots,f_r):(\C^N,0)\to(\C^r,0)$ consider a version of degeneracy
matrix,
 describing the singular locus:
\beq
N(\uf):=\bpm x_1\frac{\di f_1}{\di x_1}&\dots&x_n\frac{\di f_1}{\di x_n}&f_1&\dots&0\\
\dots&\dots&\dots&\dots&\dots&\dots\\
x_1\frac{\di f_r}{\di x_1}&\dots&x_n\frac{\di f_r}{\di x_n}&0&\dots&f_r
\epm.
\eeq
\bed
Consider the map $\uf=(f_1,\dots,f_r):(\C^N,0)\to(\C^r,0)$. Suppose $f_1,\dots,f_r$ are
convenient.
 The map $\uf$ is called \Nnd\ if
\\(i) the sequences $f_1,\dots,f_p$ are non-degenerate for any $p=1,\dots,r-1$, and
\\(ii) the submodule $N(\uf)\sset\C[[x_1,\dots,x_N]]^{\oplus r}$ is \Nnd.
\eed

\beR
 In the non-hypersurface case  this notion of non-degeneracy notion is more restrictive than the
 original definition of \cite{Kouchnirenko76}, \cite{Khovanskii1978}.
 Still, for a fixed set of diagrams it is a generic property. As we are interested only in the
 topological invariants, $\mu$, $p_g$,
 we can always assume that the complete intersection $(X,0)$ is non-degenerate in this strict
 sense, and use the formulas of the next subsection.
\eeR

\subsection{Invariants for \Nnd\ complete intersection
singularities}\label{Sec.Background.mu.pg.for.Nnd.C.I.}
\subsubsection{Milnor number for \Nnd\ singularities}
Let $(X,0)\sset(\C^N,0)$ be an isolated hypersurface singularity,
non-degenerate \wrt its \ND\ $\Ga$. Assume that $\Ga$ is convenient. In this case
the Milnor number was computed by \cite{Kouchnirenko76}:
\beq
\mu(X,0)=\suml_{0\le i\le N}(-1)^{N-i}i!\Vol_i(\Ga^-).
\eeq

The Milnor number of a \Nnd\ complete intersection singularity, $(X,0)\sset(\C^N,0)$, was
obtained in \cite[Theorem 3.9]{Bivia-Ausina2007}:
\beq\label{Eq.Milnor.Number.Nnd}\ber
\mu(X,0)=\suml^{n+r}_{j=r}(-1)^{n+r-j}\Big(\suml_{\substack{I\sseteq\{1,\dots,n+r\}\\|I|=j}}j!
a_j((\Ga^+_1)^I,\dots,(\Ga^+_r)^I)\Big)+(-1)^{n+1},\quad {\rm where}
\\
a_j((\Ga^+_1)^I,\dots,(\Ga^+_r)^I):=\suml_{\substack{\uk\in\cK_{j,r}\\k_1,\dots,k_r\ge1}}
\Covol_j\Big(((\Ga^+_1)^I)^{k_1},\dots,((\Ga^+_r)^I)^{k_r}\Big).
\eer\eeq
Here $I$ runs over all the coordinate planes, $(\Ga^+_j)^I$ is the intersection of $\Ga^+_j$
with the $I$-th coordinate plane.
 The coefficient $\Covol_j\Big(((\Ga^+_1)^I)^{k_1},\dots,((\Ga^+_r)^I)^{k_r}\Big)$ is the
 $j$-dimensional
  mixed-covolume, defined in \ref{Sec.Background.Mixed.Covolumes}.

In the particular case, when all the diagrams are proportional, i.e. $\Ga^+_j=d_j\Ga^+$,
$d_j\in\Q_{>0}$, one gets:
\beq\ber
\mu(X,0)=\suml^{n+r}_{j=r}(-1)^{n+r-j}\Big(\Theta_j(d_1,\dots,d_r)j! \Vol_j(\Ga^-)\Big)
+(-1)^{n+1},\quad {\rm where}
\\
\Theta_j(d_1,\dots,d_r):=(\prodl^r_{i=1}d_i)\suml_{\uk\in\cK_{n,r}}(\prodl^r_{i=1}d^{k_i}_i).
\eer\eeq
This was obtained in \cite[p.27]{Oka.1990}, for $\{d_i\}$ integers
and in \cite[Corollary 6.12]{Bivia-Ausina2007}, for  $\{d_i\}$ rational.
 The case when the \ND s of $f_1,\dots,f_r$ are `very close'  was clarified
  in \cite{Martin-Pfister}.
(Here `very close' means that all $f_i$ are non-degenerate with respect to the `common' diagram
$\Ga$, defined by the union of supports $\cupl_i \Supp(f_i)$.)

\subsubsection{Singularity genus for \Nnd\ singularities}
For a \Nnd\ hypersurface singularity the singularity genus is expressible as
 $p_g(X,0)=p_g(\Ga^+):=|\Ga^-\cap\Z^{n+1}_{>0}|$,
the number of strictly positive lattice point in $\Ga^-$.
 (This was proven for curves by Hodge (1928) and in higher dimensions in \cite{Merle-Teissier}, \cite{Khovanskii1978},
 \cite{Saito1981}.)

For \Nnd\ ICIS the expression for the singularity genus is the following, cf.\\
\cite{Khovanskii1978}, \cite[Theorem 2.4]{Morales1984},
\beq\label{Eq.pg.for.nnd.C.I.}
p_g(X,0)=p_g(\suml_{j=1}^r\Ga^+_j)-\suml_{i=1}^r p_g(\suml_{j\neq i}\Ga^+_j)+
\suml_{1<i_1<i_2\le r}
p_g(\suml_{j\not\in\{i_1,i_2\}}\Ga^+_j)-\cdots+(-1)^{r+1}\suml^r_{j=1}p_g(\Ga^+_j).
\eeq

\subsection{Ehrhart polynomial}\label{Sec.Background.Ehrhart.Polynomial}
Let $\De\sset\Z^N$ be a convex lattice  polytope. Let
$k\De\sset\Z^N$ be the polytope obtained by homogeneous $k$--scaling.
The number of lattice points in $k\De$ can be expressed  by the Ehrhart
polynomial of $\De$:
\beq
|k\De\cap\Z^N|=k^N \Vol_N(\De)+\frac{k^{N-1}}{2}\Vol_{N-1}(\De)+\suml_{i=1}^{N-2}c_i k^i+1,
 \eeq
 where $\Vol_N(\De)$ is the lattice $N$-dimensional
volume, $\Vol_{N-1}(\De)$ is the $(N-1)$-dimensional lattice volume of all
the top dimensional faces of $\De$. The remaining coefficients
$\{c_1,\dots,c_{N-2}\}$ are complicated. The number of lattice points
lying in the interior $\overset{\circ}{k\De}$ of $k\De$ is expressible as:
 \beq\label{Eq.Ehrhart.Pol.Polytope}
|\overset{\circ}{k\De}\cap\Z^N|=k^N
\Vol_{N}(\De)-\frac{k^{N-1}}{2}\Vol_{N-1}(\De)+\suml_{i=1}^{N-2}(-1)^{N-i}c_i k^i+(-1)^{N}.
 \eeq
 For a polygon in $\R^2$  Ehrhart formulas reduce
to the classical Pick's theorem:
 \beq\label{Eq.Ehrhart.Pol.Inner.Part.Polytope}
|k\De\cap\Z^N|=k^2\Vol_2(\De)+\frac{k}{2}\Vol_{1}(\De)+1,\quad\quad\quad
|\overset{\circ}{k\De}\cap\Z^N|=k^2\Vol_2(\De)-\frac{k}{2}\Vol_{1}(\De)+1.
\eeq

\bex\label{Ex.Computing.pg(Gamma).general}
To obtain the expression for $p_g(\Ga)=p_g(\Ga^-)$, i.e. the number of $\Z^N_{>0}$ points on or
under $\Ga$,
present $\Ga^-=\De\smin(\Ga^+\cap\De)$. Here $\De$ is a large enough convex polytope that lies
in $\R^N_{\ge0}$
 and contains $\Ga^-$. Then:
 $p_g(\Ga)=|\overset{\circ}{\De}\cap\Z^N|-|\overset{\circ}{\overline{\Ga^+\cap\De}}\cap\Z^N|$.
 Equation (\ref{Eq.Ehrhart.Pol.Polytope}) gives:
\beq\label{Eq.pg.Ehrhart.expansion}
p_g(\Ga)=p_g(\Ga^-)=\Vol_N(\Ga^-)+\frac{\Vol_{N-1}(\Ga)-\Vol_{N-1}(\Ga^-)}{2}+\suml^{N-2}_{i=1}
(-1)^{N-i}k^i\Big(c_i(\De)-c_i(\Ga^+\cap\De)\Big).
\eeq
\eex

\section{Proof of the bound for large enough \Nnd\ complete intersections}
Here we prove theorem \ref{Thm.Bound.for.Nnd.C.I.}.
The proof goes in 2 steps. First, we reduce the problem to a combinatorial statement, by
 expressing $\mu$ and $p_g$ in terms of the (mixed-)covolumes,
  $\{\Covol((\Ga^+_1)^{k_1},\dots,(\Ga^+_r)^{k_r})\}_{(k_1,\dots,k_r)}$.
 Then we compare the leading terms of $\mu$ and $p_g$ and prove $\lt(\mu)\ge \lt(p_g)$ with
 equality
 only in the case $\Ga_1=\cdots=\Ga_r$. This proves the theorem when not all the diagrams
 coincide.
 Finally, in the case $\Ga_1=\cdots=\Ga_r$, we prove the theorem by comparison of the
 second-order terms.

{\bf Step 1.} Consider the isolated complete
 intersection singularity, $(X,0)=\{f_1=\cdots=f_r=0\}\sset(\C^{n+r},0)$, \Nnd\ \wrt\ the
 diagrams $(\Ga_1,\dots,\Ga_r)$.
 Suppose all the diagrams are convenient. The expressions for $\mu(X,0)$, $p_g(X,0)$ are
 given in \S\ref{Sec.Background.mu.pg.for.Nnd.C.I.}.

 We assume all the diagrams $\Ga_i$ to be large enough,
in particular $\Vol_{n+r}(\Ga^-_i)\gg \Vol_{n+r-1}(\Ga^-_i)\gg\cdots$.
 Thus in the comparison of $\mu$ vs $p_g$ it is enough to compare only the higher order terms.
First we compare the leading terms.

The leading term for $p_g(\Ga)$ is obtained from Ehrhart expansion, equation
(\ref{Eq.pg.Ehrhart.expansion}):
 $\lt(p_g(\Ga))=\Vol_{n+r}(\Ga^-)$. Thus equation (\ref{Eq.pg.for.nnd.C.I.}) gives:
\beqm
\lt(p_g(X,0))=
\Covol_{n+r}(\suml_{j=1}^r\Ga^+_j)-\suml_{i=1}^r \Covol_{n+r}(\suml_{j\neq i}\Ga^+_j)+\\+
\suml_{1<i_1<i_2\le r}
\Covol_{n+r}(\suml_{j\not\in\{i_1,i_2\}}\Ga^+_j)-\cdots+(-1)^{k+1}\suml^r_{j=1}\Covol_{n+r}(\Ga^+_j).
\end{multline}
Expand all the brackets using mixed covolumes, \S\ref{Sec.Background.Mixed.Covolumes}, to get:
\beqm
\lt(p_g(X,0))=
\suml_{\uk\in\cK_{n+r,r}}\bin{n+r}{k_1,\dots,k_r}\Covol_{n+r}\Big((\Ga^+_1)^{k_1},\dots,(\Ga^+_r)^{k_r}\Big)-
\\-\suml_{i=1}^r \suml_{\substack{\uk\in\cK_{n+r,r}\\k_i=0}}\bin{n+r}{k_1,\dots,k_r}
\Covol_{n+r}\Big((\Ga^+_1)^{k_1},\dots,(\Ga^+_r)^{k_r}\Big)
+\cdots +(-1)^{k+1}\suml^r_{j=1}\Covol_{n+r}(\Ga^+_j).
\end{multline}
By the exclusion-inclusion formula, equation (\ref{Eq.Exclusion.Inclusion.Formula}), we get:
\beq
\lt(p_g(X,0))=\suml_{\substack{\uk\in\cK_{n+r,r}\\k_1,\dots,k_r\ge1}}\bin{n+r}{k_1,\dots,k_r}
\Covol_{n+r}\Big((\Ga^+_1)^{k_1},\dots,(\Ga^+_r)^{k_r}\Big).
\eeq
The leading term of $\mu$ is immediate:
\beq
\lt(\mu(X,0))=(n+r)!\suml_{\substack{\uk\in\cK_{n+r,r}\\k_1,\dots,k_r\ge1}}
\Covol\Big((\Ga^I_1)^{k_1},\dots,(\Ga^I_r)^{k_r}\Big).
\eeq
To prove the initial equality it is enough to check $\lt(\mu(X,0))>\cf\cdot \lt(p_g(X,0))$.
 We prove:
\beq
\lt(\mu(X,0))\ge\cf\cdot\lt(p_g(X,0)), \text{ and equality occurs iff $\Ga_1=\cdots=\Ga_r$}.
\eeq
 (For example, the equality occurs in the hypersurface case, $r=1$.)
But this is exactly the inequality presented in equation (\ref{Eq.Ineq.of.mixed.Covolumes}),
proved in \cite[\S4]{Kerner-Nemethi.Combinatorial.Lemma}.

{\bf Step 2.} The comparison of the leading terms, as above, proves $\mu> \cf p_g$ when at least
two
 diagrams among $\{\Ga_i\}$ do not coincide. It remains to check the case $\Ga_1=\cdots=\Ga_r$.
In this case the expressions for $\mu$ and $p_g$ simplify:
\beq\ber
\mu(X,0)=\suml^{n+r}_{j=r}(-1)^{n+r-j}|\cK_{j-r,r}|j! \Vol_j(\Ga^-)+(-1)^{n+1},\\\\
p_g(X,0)=p_g(r\Ga^+)-r p_g((r-1)\Ga^+)+\bin{r}{2}p_g((r-2)\Ga^+)-\cdots+(-1)^{r+1}rp_g(\Ga^+)
\eer\eeq
Now the expansions by the orders of $\Ga$ are:
\beq\ber
\mu(X,0)=(n+r)!|\cK_{n,r}|\Vol_{n+r}(\Ga^-)-(n+r-1)!|\cK_{n-1,r}|\Vol_{n+r-1}(\Ga^-)+\cdots,
\\\\
p_g(X,0)=\suml^r_{j=0}(-1)^j\bin{r}{j}\Big(\Vol_{n+r}((r-j)\Ga^-)+\frac{\Vol_{n+r-1}((r-j)\Ga)-\Vol_{n+r-1}((r-j)\Ga^-)}{2}
\Big)+\cdots
\eer\eeq
Note that $\Vol_i((r-j)\Ga^-)=(r-j)^i \Vol_i(\Ga^-)$ and
$\Vol_{n+r-1}((r-j)\Ga)=(r-j)^{n+r-1}\Vol_{n+r-1}(\Ga)$.
 Thus one has:
\beq
p_g(X,0)=\stir{n+r}{r}\Vol_{n+r}(\Ga^-)+\stir{n+r-1}{r}\frac{\Vol_{n+r-1}(\Ga)-\Vol_{n+r-1}(\Ga^-)}{2}+\cdots
\eeq
Thus we need  to prove:
\beqm
(n+r)!|\cK_{n,r}|\Vol_{n+r}(\Ga^-)-(n+r-1)!|\cK_{n-1,r}|\Vol_{n+r-1}(\Ga^-)>\\\hspace{3cm}>\cf\Big(
\stir{n+r}{r}\Vol_{n+r}(\Ga^-)+\stir{n+r-1}{r}\frac{\Vol_{n+r-1}(\Ga)-\Vol_{n+r-1}(\Ga^-)}{2}\Big).
\end{multline}
The leading terms here cancel. (This was shown in Step 1 and can be also checked explicitly:
 $(n+r)!|\cK_{n,r}|\Vol_{n+r}(\Ga)=\cf\stir{n+r}{r}\Vol_{n+r}(\Ga^-)$.)
  Therefore it remains to prove:
\beq
-(n+r-1)!|\cK_{n-1,r}|\Vol_{n+r-1}(\Ga^-)>\cf\stir{n+r-1}{r}\frac{\Vol_{n+r-1}(\Ga)-\Vol_{n+r-1}
(\Ga^-)}{2}.
\eeq
Use the definition of $\cf$ to present this in the form:
$$\frac{\Vol_{n+r-1}(\Ga^-)-\Vol_{n+r-1}(\Ga)}{2}\frac{(n+r)(n+r-1)}{n}
\frac{\stir{n+r-1}{r}}{\stir{n+r}{r}}>\Vol_{n+r-1}(\Ga^-).$$
We claim that $\Vol_{n+r-1}(\Ga)\le\frac{\Vol_{n+r-1}(\Ga^-)}{(n+r)}$. This can be seen, e.g. by
the projection of $\Ga$ on all
the coordinate hyperplanes, $\{x_i=0\}^{n+r}_{i=1}$.
Substitute this inequality and cancel $\Vol_{n+r-1}(\dots)$ It remains to prove:
$$\stir{n+r-1}{r}\Big/\stir{n+r}{r}>\frac{2n}{(n+r-1)^2}.$$
For $r=1$, $n>2$ this inequality is verified directly: $1>\frac{2}{n}$.
 Thus we assume $r>1$ and use the recurrence relations of
 \S\ref{Sec.Stirling.Coefficients.2'nd.kind}.
 This gives:
\beq
\frac{\stir{n+r-1}{r}}{\stir{n+r}{r}}=\frac{1}{r+\frac{\stir{n+r-2}{r-1}}{\stir{n+r-1}{r}}}>\frac{1}{r+1}
\eeq
(For the later inequality see \S\ref{Sec.Stirling.Coefficients.2'nd.kind}.)
Therefore, it is enough to check: $\frac{1}{r+1}>\frac{2n}{(n+r-1)^2}$.
Note that $\frac{1}{r+1}-\frac{2n}{(n+r-1)^2}=\frac{(r-1)^2+(n^2-4n)}{(r+1)(n+r-1)^2}$.
This leaves only one case to check separately: $(r,n)=(2,3)$. In this case:
\beq
\frac{\stir{n+r-1}{r}}{\stir{n+r}{r}}=\frac{\stir{3+2-1}{2}}{\stir{3+2}{2}}=\frac{7}{15}>
\frac{6}{16}=\frac{2\cdot 3}{(3+2-1)^2}=\frac{2n}{(n+r-1)^2}.\quad\bull
\eeq

\section{A stronger asymptotic bound for hypersurfaces}\label{Sec.Stronger.Bound.Hypersurfaces}
The proof of theorem \ref{Thm.Asympt.Bound.Hypersurfaces} is in \S\ref{Sec.Asympt.Bound.Hypersurfaces}.
Although the germ $(X,0)\sset(\C^{n+1},0)$ is a local object, the statement of the theorem contains the projective hypersurface $\P T_{(X,0)}\sset\P^n$. In \S\ref{Sec.Khovanskii.Kouchnirenko.Formula.Proj.Hypers} we derive some facts about the Milnor number $\mu(\P T_{(X,0)})$.

\subsection{An auxiliary  Khovanskii-Kouchnirenko type formula.}\label{Sec.Khovanskii.Kouchnirenko.Formula.Proj.Hypers}
 Let $\De\sset\R^{n+1}_{\ge0}$ be a convex  lattice
polytope such that $\Span_\R(\De)=\R^{n+1}$. Let $(\C^*)^{n+1}\sset
Y_\De$ be the corresponding toric completion, with the natural sheaf
$\cO_{Y_\De}(1)$. Let $D_\infty:=Y_\De\smin(\C^*)^{n+1}$ be the
divisor at infinity. The variety $Y_\De$ is in general non-smooth.
Suppose it is smoothable, i.e. there exists a flat family
$(\cY,\cL_\cY)$ over $(\C^1,0)$ such that
$\cY|_0=Y_\De$, $\cY|_{t\neq0}$ is smooth and $\cL|_{\pi^{-1}(0)}=\cO_{Y_\De}(1)$.

Let $\De_0\sset\De$ be a lattice sub-polytope, let $f$ be a
function supported on $\De_0$ and non-degenerate \wrt $\De_0$. Let
$X_{\De_0}=\overline{\{f=0\}}\sset Y_\De$ be the corresponding
projective hypersurface. By construction all its singularities lie
on the boundary $D_\infty$. Note that $X_{\De_0}$ can have non-isolated singularities.

Suppose $\dim(\De_0)=n+1$, in particular $\Span_\R(\De_0)=\R^{n+1}$.
Let $X_\De$ be a generic (partial) smoothing of $X_{\De_0}$  inside $Y_\De$. Namely, $X_\De\sset
Y_\De$ is
 a hypersurface, defined by $\{f_t=0\}$, where $\Supp(f_t)=\De$ and $f_t$ is non-degenerate on
 $\De$.
 So $X_\De\cap(\C^*)^{n+1}$ is smooth and $X_\De$ intersects $D_\infty$ transversally.
 Note that $X_\De$ itself is smoothable,  in the  family $(\cY,\cL_\cY)$, and its smoothing is
also a smoothing of $X_{\De_0}$.

Define the Milnor number,  $\mu(X_{\De_0}):=(-1)^n\Big(\chi(X_\De)-\chi(X_{\De_0})\Big)$.

If $Y_\De$ is itself smooth then $X_\De$ is smooth and this definition coincide with that of
equation (\ref{Eq.Def.of.Mu.Pg}).

\bel\label{Thm.Bound.Total.Milnor.Number.Area.Polytope} Under the
assumptions as above:
$\mu(X_{\De_0})=(n+1)!\Vol_{n+1}(\De\smin\De_0)-\mu\Big(X_{\De_0}\cap
D_\infty\Big)$.
 \eel
\bpr By
\cite[pg. 59]{Khovanskii1978}
\beq
\chi(X_\De\cap(\C^*)^{n+1})=(-1)^{n}(n+1)!\Vol_{n+1}(\De)
\eeq
and similarly for $X_{\De_0}$. (Here we use the assumption that $\dim(\De_0)=n+1$.)
Hence
\beq
\mu(X_{\De_0})=(n+1)!\Vol_{n+1}(\De\smin\De_0)+
(-1)^n\Big(\chi(X_\De\cap D_\infty)-\chi(X_{\De_0}\cap
D_\infty)\Big) \eeq Finally, as $X_\De$ intersects $D_\infty$
transversally and $X_\De\cap(\C^*)^{n+1}$ is smooth we obtain that
if $X_\ep$ is a smoothing of $X_{\De_0}$ then
$\chi(X_\ep\cap(\C^*)^{n+1})=\chi(X_\De\cap(\C^*)^{n+1})$, and
 $X_\ep\smin(X_\ep\cap(\C^*)^{n+1})$ is a smoothing of $X_\De\cap D_\infty$.
Thus
\beq
\mu(X_\De)=(-1)^n\Big(\chi(X_\ep\smin(X_\ep\cap(\C^*)^{n+1})-\chi(X_\De\cap
D_\infty)\Big),
\eeq
 and
\beq
\mu(X_{\De_0})=(-1)^n\Big(\chi(X_\ep\smin(X_\ep\cap(\C^*)^{n+1})-\chi(X_{\De_0}\cap
D_\infty)\Big).
\eeq
\epr
\bex In the simplest case, suppose $\De=\Conv(x^p_0,\dots,x^p_{n+1})\sset\R^{n+2}$, so that
$(Y_\De,\cO_{Y_\De}(1)\approx(\P^{n+1},\cO_{\P^{n+1}}(p))$. Suppose $\De_0$ intersects all the
(one-dimensional) edges of $\De$, then $X_{\De_0}$ has only isolated
singularities. Then iterating the formula of the lemma we get
Kouchnirenko's formula:
\beq
\mu(X_{\De_0})=(n+1)!\Vol_{n+1}(\De\smin\De_0)-(n)!\Vol_{n}(\De\smin\De_0)+\dots.
\eeq
This formula will be used in equation (\ref{Eq.using.Milnor.number.Proj.Tang.Cone}).
\eex

\subsection{Proof of theorem \ref{Thm.Asympt.Bound.Hypersurfaces}.}\label{Sec.Asympt.Bound.Hypersurfaces}

Let $(X,0)\sset(\C^{n+1},0)$ be an isolated hypersurface singularity, non-degenerate \wrt its diagram $\Ga_{(X,0)}$.

By direct check, if $(X,0)$ is a homogeneous isolated hypersurface singularity (and thus $\mu(\P T_{(X,0)})=0$), we have the equality:
\beq
\mu(X,0)-\Big((p-1)^{n+1}-\frac{p!}{(p-n-1)!}\Big)=(n+1)!p_g(X,0).
\eeq
Therefore we assume that $(X,0)$ is not an ordinary multiple point, in particular $\P T_{(X,0)}$ is not smooth.

The combinatorial formulas for Milnor number and
geometric genus of a \Nnd\ singularity are given in \S\ref{Sec.Background.mu.pg.for.Nnd.C.I.}.
We want to prove: for any
\ND\ $\Ga$ there exists $k_0$ such that for $k\ge k_0$ one has
\beq\label{Eq.Ineq.To.be.proved.First}
 \mu(k\Ga)-\mu(\P T_{(kX,0)})-\Big((kp-1)^{n+1}-\frac{(kp)!}{(kp-n-1)!}\Big)>
(n+1)!p_g(k\Ga^-).
\eeq
(Here $kX$ denotes the corresponding projective hypersurface.
If the singularities of $\P T_{(kX,0)}$ are non-isolated then the term $\mu(\P T_{(kX,0)})$ is omitted.)
As in the proof of theorem \ref{Thm.Bound.for.Nnd.C.I.} we expand the whole expression in
powers of $k$ and prove that the leading term is positive.

{\bf Step 1.} Equation (\ref{Eq.pg.Ehrhart.expansion}) gives:
\beq\label{Eq.pg.expansion.Ehrhart}
p_g(k\Ga^-)=k^{n+1}\Vol_{n+1}(\Ga^-)+\frac{k^{n}}{2}\Big(\Vol_{n}\Ga-\Vol_{n}\Ga^-\Big)+\lot
\eeq
The Kouchnirenko formula for Milnor number gives:
\beq
\mu(k\Ga)=k^{n+1}(n+1)!\Vol_{n+1}(\Ga^-)-k^{n}n!\Vol_{n}(\Ga^-)+\lot
\eeq
If the singularities of $\P T_{(X,0)}$ are isolated then in particular $\dim(\De_0)=n$.
Then lemma \ref{Thm.Bound.Total.Milnor.Number.Area.Polytope} reads:
\beq\label{Eq.using.Milnor.number.Proj.Tang.Cone}
\mu(\P T_{(kX,0)})=n!\Vol_n(k\De\smin\De_0)-\mu(X_{k\De_0}\cap D_\infty)\quad
 \text{ and  the $k$-order of }\mu(X_{k\De_0}\cap D_\infty) \text{ is lower than $n$.}
\eeq
Here $\De=\Conv(x^p_1,\dots,x^p_{n+1})$, while $\De_0$ is the Newton polyhedron of $\P T_{(X,0)}$. In what follows we denote $\De_0$ by $\Ga(\P T_{(X,0)})$.

Finally, expand \beq
(kp-1)^{n+1}-\frac{(kp)!}{(kp-n-1)!}=\frac{(n+1)(n-2)}{2}p^{n}k^{n}-\bin{n+2}{3}\frac{3n-7}{4}p^{n-1}k^{n-1}+\lot
\eeq
 Substitute all the data into the inequality
(\ref{Eq.Ineq.To.be.proved.First}) to get the expansion:
 \beq
\frac{k^{n}(n+1)!}{2}\Bigg(\frac{n-1}{n+1}\Vol_{n}\Ga^--\Vol_{n}\Ga-\frac{n-2}{n!}p^{n}-
\frac{2}{n+1}\Vol_{n}\Big(\De\smin\Ga(\P
T_{(X,0)})\Big)\Bigg)+\lot
\eeq

To prove that this expression is positive/non-negative we check the
coefficient of $k^{n}$. If the singularities of $\P T_{(X,0)}$ are non-isolated then we can omit the term
$\De\smin\Ga(\P T_{(X,0)})$. However we prove the non-negativity even with that term. (Note that $\Vol_{n}\Big(\De\smin\Ga(\P T_{(X,0)})\Big)$ is non-negative.)

Since $\frac{p^n}{n!}$ is the volume $\Vol_{n}\De$, we need to prove
\beq\label{Eq.Ineq.To.be.proved.Second}
\frac{n-1}{n+1}\Vol_{n}\Ga^--\Vol_{n}\Ga+\frac{2}{n+1}\Vol_{n}\Ga(\P
T_{(X,0)})-(n-2+\frac{2}{n+1})\Vol_{n}\De>0.
\eeq

{\bf Step 2.} Let
$\Ga=\cupl_\al\si_\al$ be the decomposition into the top-dimensional
faces. Here $\al$ belongs to some set and we fix a special value
$\al=p$ by $\si_p:=\Ga\cap \De$. If $\si_p$  is not
top-dimensional, then it is omitted.

Let $\pi_j: \R^{n+1}\to\{x_j=0\}\sset\R^{n+1}$ be the projection
onto a coordinate hyperplane. Note that $\pi_j$ sends $\Z^{n+1}$ to
$\Z^{n}$, in particular $\pi_j(\si_\al)$ is a lattice polytope.
Consider the union of the images of such projections,
$\pi\si_\al=\cupl_j\pi_j(\si_\al)$.

Now, we return  to inequality (\ref{Eq.Ineq.To.be.proved.Second}). We have
 \beq
\Vol_{n}\Ga=\suml_{\al\neq p}\Vol_{n}\si_\al+\Vol_{n}\si_p,\quad\quad
\Vol_{n}\Ga^-=\suml_{\al\neq p}\Vol_{n}\pi\si_\al+\Vol_{n}\pi\si_p.
\eeq
Here the sums $\suml_{\al\neq p}(\dots)$ are non-empty as $(X,0)$ is not an ordinary multiple point.

Note that $\Vol_{n}\pi\si_p=(n+1)\Vol_{n}\si_p$ and $\Ga(\P
T_{(X,0)})=\si_p$.
Thus the inequality (to be proved) becomes:
\beq
\suml_{\al\neq
p}\Big(\frac{n-1}{n+1}\Vol_{n}\pi\si_\al-\Vol_{n}\si_\al\Big)-
(n-2+\frac{2}{n+1})\Vol_{n}(\De\smin\si_p)>0.
\eeq

{\bf Step 3.}
Consider the projection $\Ga\norm\De$ defined by $pt\to \De\cap
 \rm{line}(0,pt)$. This projection is surjective as a map of points of
$\Ga$ with {\em real} coordinates. In general the lattice points of
$\Ga$ are not sent to the lattice points of $\De$.

The image of a face, $\nu(\si_\al)\sset\De$ is a {\em rational}
polytope. Let $\Vol_{n}^\R(\nu(\si_\al))$ denote its rational {\em
normalized} volume, namely:
$\Vol_{n}^\R(\nu(\si_\al)):=\Vol_{n}^\R(\pi_j\nu(\si_\al))$, for any
$j$. Here $\Vol_{n}^\R(\pi_j\nu(\si_\al))$ is the usual volume in the
hyperplane $\R^{n-1}$. Note that
$\Vol_{n}(\De\smin\si_p)=\suml_{\al\neq p}\Vol_{n}^\R(\nu(\si_\al))$.
Thus the inequality can be written in the form \beq \suml_{\al\neq
p}\Big(\frac{n-1}{n+1}\Vol_{n}\pi\si_\al-\Vol_{n}\si_\al-
(n-2+\frac{2}{n+1})\Vol_{n}^\R(\nu(\si_\al))\Big)>0. \eeq We prove
that each summand is positive.

{\bf Step 4.} Suppose that the top dimensional face $\si_\al$ lies in the hyperplane
$\suml_{j=1}^{n+1}\frac{x_j}{a_j}=\const$.
 Here $\{a_j\}$ are natural numbers and $\gcd(a_1,\dots,a_{n+1})=1$. Then the primitive normal to
 the face has coordinates:
 $\cN_\al=(\frac{\prod a_i}{a_1d},\dots,\frac{\prod a_i}{a_{n+1}d})$, where
 $d:=\gcd(\frac{\prod a_i}{a_1},\dots,\frac{\prod a_i}{a_{n+1}})$.

Note that $\frac{\Vol_{n}\pi_j\si_\al}{\Vol_{n}\si_\al}=\frac{\prod
a_i}{a_jd}$. This can be obtained by comparing the lattice areas of
the simplex $\Conv(x^{a_1}_1,\dots,x^{a_{n+1}}_{n+1})$  and its
projections. Therefore \beq
\frac{n-1}{n+1}\Vol_{n}\pi\si_\al-\Vol_{n}\si_\al=\Vol_{n}\si_\al
\Big(\frac{n-1}{n+1}\suml_j\frac{\prod
a_i}{a_jd}-1\Big).
\eeq
Now compare $\Vol_{n}\si_\al$ to
$\Vol_{n}^\R\nu\si_\al$. We claim $\Vol_{n}^\R\nu\si_\al<\min_j
\Vol_{n}\pi_j\si_\al$ (note that the inequality is strict). Indeed,
the left hand side was defined (in Step 3.) as the real area
$\Vol^\R_{n-1}\pi_j\nu\si_\al$. But
$\Vol^\R_{n-1}\pi_j\nu\si_\al<\Vol_{n}\pi_j\si_\al$.

Thus $\Vol^\R_{n-1}\nu\si_\al<(\min_j \frac{\prod a_i}{a_jd})
\Vol_{n}\si_\al$. Therefore it is enough to prove the following
arithmetic statement, for $(a_1,\dots,a_{n+1})\neq(1,\dots,1)$: \beq
\frac{n-1}{n+1}\suml_j\frac{\prod
a_i}{a_jd}-1-(n-2+\frac{2}{n+1})\min_j\frac{\prod a_i}{a_jd}\ge0.
\eeq
Note that now the inequality to be proved is {\em non-strict}.
Present it in the form:
\beq
\frac{n-1}{n+1}\Big(\suml^{n+1}_{j=1} \frac{\prod a_i}{a_jd}-n\cdot
\min\limits_j \frac{\prod a_i}{a_jd}\Big)\ge1.
\eeq
As $\cN\neq(1,\dots,1)$ we have:
$\suml^{n+1}_{j=1} \frac{\prod a_i}{a_jd}\ge (n+1)\cdot\min\limits_j \frac{\prod a_i}{a_jd}+1$.
So, the inequality becomes
\beq
\frac{n-1}{n+1}\Big(1+\min\limits_j \frac{\prod
a_i}{a_jd}\Big)\ge1,
\eeq
 which is obvious for $n\ge3$.
(Just note: $1+\min\limits_j \frac{\prod a_i}{a_jd}\ge2$ and $\frac{n-1}{n+1}\cdot2\ge1$.)
\epr


\begin{thebibliography}{99}
\bibitem[Abramowitz-Stegun]{Abramowitz-Stegun}
M. Abramowitz, I.A. Stegun, {\em Handbook of mathematical functions with formulas, graphs, and
mathematical tables.}
 National Bureau of Standards Applied Mathematics Series, 55, 1964.


\bibitem[AGLV]{AGLV} V.I. Arnol'd, V.V. Goryunov, O.V. Lyashko, V.A. Vasil'ev, {\em
Singularity theory.I.} Reprint of the original English edition
from the series Encyclopaedia of Mathematical Sciences [Dynamical
systems. VI, Encyclopaedia Math. Sci., 6, Springer, Berlin, 1993].
Springer-Verlag, Berlin, 1998.

\bibitem[Artin1966]{Artin-66} M. Artin, {\it On isolated rational singularities of surfaces}.
 Amer. J. Math. 88 1966 129--136.

 \bibitem[Ashikaga1992]{Ashikaga92} T. Ashikaga, {\em Normal two-dimensional hypersurface triple
     points
 and the Horikawa type resolution.} Tohoku Math. J. (2) 44 (1992), no. 2, 177--200.

\bibitem[Berline-Vergne-2007]{Berline-Vergne-2007} N. Berline, M. Vergne, {\em Local
    Euler-Maclaurin formula
for polytopes.} Mosc. Math. J. 7 (2007), no. 3, 355--386, 573.


\bibitem[Bivi\`{a}-Ausina2004]{Bivia-Ausina2004} C. Bivi\`{a}-Ausina, {\em The integral closure
    of modules,
Buchsbaum-Rim multiplicities and Newton polyhedra.} J. London Math. Soc. (2) 69 (2004), no. 2,
407--427.

\bibitem[Bivi\`{a}-Ausina2007]{Bivia-Ausina2007} C. Bivi\`{a}-Ausina, {\em Mixed Newton numbers
    and isolated
 complete intersection singularities.} Proc. Lond. Math. Soc. (3) 94 (2007), no. 3, 749--771.


\bibitem[Buchweitz-Greuel1980]{Buchweitz-Greuel-1980}R.-O. Buchweitz, G.-M. Greuel, {\em
The Milnor number and deformations of complex curve singularities.} Invent. Math. 58 (1980), no.
3, 241--281.


\bibitem[Dimca]{Dimca92} A. Dimca, {\it Singularities and topology of hypersurfaces.}
    Universitext. Springer-Verlag,  New York, 1992.

 \bibitem[Damon1989]{Damon1989} J. Damon, {\em Topological invariants of $\mu$-constant
     deformations of complete intersection singularities.}  Quart. J. Math. Oxford Ser. (2) 40
     (1989), no. 158, 139--159.

\bibitem[Durfee1978]{Durfee1978} A.H. Durfee, {\it The signature of smoothings of complex
    surface singularities.}
Math. Ann. 232 (1978), no. 1, 85--98.


 \bibitem[Greuel1975]{Greuel1975} G.M. Greuel, {\em Der Gauss-Manin-Zusammenhang isolierter
     Singularit\"{a}ten von
 vollst\"{a}ndigen Durchschnitten.} Math. Ann. 214 (1975), 235-266.

 \bibitem[Greuel-Hamm1978]{Greuel-Hamm1978} G.M. Greuel, H.A. Hamm, {\em Invarianten
     quasihomogener vollst\"{a}ndiger
 Durchschnitte.}, Invent. Math. 49 (1978), no. 1, 67--86.

 \bibitem[Hamm1986]{Hamm1986} H.A. Hamm, {\em Invariants of weighted homogeneous singularities.}
     Journ\'{e}es Complexes 85 (Nancy, 1985), 613, Inst. \'{E}lie Cartan, 10, Univ. Nancy,
     Nancy, 1986.

\bibitem[Hamm2011]{Hamm2011} H.A. Hamm, {\em Differential forms and Hodge numbers for toric
    complete intersections},  arXiv:1106.1826.


\bibitem[Jordan1965]{Jordan1965} Ch. Jordan, {\em Calculus of finite differences}. Third
    Edition.  Introduction by Harry C. Carver Chelsea Publishing Co., New York 1965.

\bibitem[Kantor-Khovanskii-1993]{Kantor-Khovanskii1993}
 J.M. Kantor, A. Khovanskii, {\em Une application du th\'{e}or\`{e}me de Riemann-Roch
 combinatoire au polyn\^{o}me
d'Ehrhart des polytopes entiers de ${\R}^d$} C. R. Acad. Sci. Paris Sr. I Math. 317 (1993),
no. 5, 501--507.

\bibitem[Kaveh-Khovanskii-2013-1]{Kaveh-Khovanskii-Conv.Bod}  K. Kaveh, A.G. Khovanskii, {\em
    Convex bodies and multiplicities of ideals}, arXiv:1302.2676.

 \bibitem[Kaveh-Khovanskii-2013-2]{Kaveh-Khovanskii-Mix.Mult}  K. Kaveh, A.G. Khovanskii, {\em
     On mixed multiplicities of ideals}, arXiv:1310.7979.


\bibitem[Kerner-N\'emethi2009]{Kerner-Nemethi2009} D. Kerner and A. N\'emethi, {\em The Milnor
    fibre signature is not semi-continuous},
 Proc. of the Conference ``Topology of Algebraic Varieties'', Jaca (Spain),  Contemporary Math.
 538 (2011), 369--376.



\bibitem[Kerner-N\'emethi2011]{Kerner-Nemethi.Durfee1} D. Kerner, A. N\'emethi,
{\em A counterexample to Durfee conjecture}, Comptes Rendus Math\'ematiques
 de l'Acad\'emie des Sciences, vol.34 (2012), no.2.  arXiv:1109.4869.

\bibitem[Kerner-N\'emethi2013]{Kerner-Nemethi.Durfee2} D. Kerner, A. N\'emethi,
{\em The 'corrected Durfee's inequality' for homogeneous complete intersections},
 Mathematische Zeitschrift, Math. Z. 274 (2013), no. 3--4, pp.1385--1400.

\bibitem[Kerner-N\'emethi2014]{Kerner-Nemethi.Combinatorial.Lemma} D. Kerner, A. N\'emethi,
{\em A generalized FKG--inequality for compositions}, to appear in Journal of Combinatorial Theory, Series A.
 arXiv: 1412.8200.




\bibitem[Khovanskii1978]{Khovanskii1978} A.G. Khovanskii,
{\em Newton polyhedra, and the genus of complete intersections.}
 (Russian) Funktsional. Anal. i Prilozhen. 12 (1978), no. 1, 51--61.


\bibitem[Kouchnirenko1976]{Kouchnirenko76} A.G. Kouchnirenko, {\it Poly\`{e}dres de Newton et
    nombres de Milnor.}
Invent. Math. 32 (1976), no. 1, 1--31.


\bibitem[Laufer1977]{Laufer1977} H.B. Laufer, {\em  On $\mu$ for surface singularities,}
 Proceedings of Symposia in Pure Math.  30, 45--49,  1977.

\bibitem[Looijenga]{Looijenga-book} E. Looijenga, {\em Isolated Singular Points on Complete
    Intersections.}
 London Math. Soc. LNS 77, CUP, 1984.

\bibitem[Looijenga1986]{Looijenga1986} E. Looijenga, {\em Riemann-Roch and smoothings of
    singularities.}
 Topology 25 (1986), no. 3, 293--302.


\bibitem[Martin-Pfister]{Martin-Pfister} B. Martin, G. Pfister, {\em Milnor number of complete
    intersections and Newton polygons.} Math. Nachr. 110 (1983), 159--177.


\bibitem[Melle-Hern\'{a}ndez2000]{Melle00} A. Melle-Hern\'{a}ndez, {\it Milnor numbers for
    surface singularities.}
  Israel J. Math.  115  (2000), 29--50.

\bibitem[Merle-Teissier]{Merle-Teissier} M. Merle,  B. Teissier {\em Conditions d'adjonction, d'apr\`{e}s DuVal} in S\'{e}minaire sur les Singularit\'{e}s des Surfaces. Lecture Notes in Mathematics, 777. Springer, Berlin, 1980. 229-245.

\bibitem[Milnor-book]{Milnor-book} J. Milnor, {\em Singular points of complex hypersurfaces,}
 Annals of Math. Studies  61, Princeton University Press 1968.

\bibitem[Morales1984]{Morales1984}  M. Morales, {\em Poly\`{e}dre de Newton et genre
    g\'{e}om\'{e}trique d'une singularit\'{e} intersection compl\`{e}te.}
 Bull. Soc. Math. France  112  (1984),  no. 3, 325--341.


 \bibitem[Morales1985]{Morales1985} M. Morales, {\em Fonctions de Hilbert, genre
     g\'{e}om\'{e}trique
d'une singularit\'{e} quasi-homog\`{e}ne Cohen-Macaulay.} C. R. Acad. Sci. Paris S\'{e}r. I
Math. 301 (1985),
no. 14, 699--702.



\bibitem[N\'emethi98]{Nemethi98} A. N\'emethi, {\em Dedekind sums and the signature of
    $f(x,y)+z^N$}.
Selecta Math. (N.S.) 4 (1998), no. 2, 361--376.

\bibitem[N\'emethi99]{Nemethi99} A. N\'emethi,  {\em  Dedekind sums and the signature of
    $f(x,y)+z^N$,II.}
 Selecta Math. (N.S.)  5 (1999), 161--179.


\bibitem[Oka.1990]{Oka.1990} M. Oka, {\em Principal zeta-function of nondegenerate complete
    intersection singularity.} J. Fac. Sci. Univ. Tokyo Sect. IA Math. 37 (1990), no. 1, 11--32.



\bibitem[Oka]{Oka-book} M. Oka, {\em Non-degenerate complete intersection singularity.}
    Actualit\'{e}s Math\'{e}matiques. [Current Mathematical Topics] Hermann, Paris, 1997.


\bibitem[Saito1981]{Saito1981} M. Saito, {\em On the exponents and the geometric genus of an
    isolated hypersurface
 singularity}. Singularities, Part 2 (Arcata, Calif., 1981), 465--472,
 Proc. Sympos. Pure Math., 40, Amer. Math. Soc., Providence, RI, 1983.

\bibitem[Seade-book]{Seade} J. Seade, {\it On the Topology of Isolated Singularities in Analytic
    Spaces.} Progress in
Mathematics 241, Birkh\"auser 2006.

\bibitem[Teissier1978]{Teissier1978} B. Teissier, {\em On a Minkowski-type inequality for
    multiplicities. II.}
 C. P. Ramanujam a tribute, pp. 347--361, Tata Inst. Fund. Res. Studies in Math., 8, Springer,
 Berlin-New York, 1978.

\bibitem[Teissier2004]{Teissier2004} B. Teissier, {\em Monomial ideals, binomial ideals,
    polynomial ideals.}
  Trends in commutative algebra, 211--246, Math. Sci. Res. Inst. Publ., 51, Cambridge Univ.
  Press, Cambridge, 2004.


\bibitem[Tomari1993]{Tomari93} M. Tomari, {\em The inequality $8p_g<\mu$ for hypersurface
    two-dimensional
 isolated double points.} Math. Nachr. 164 (1993), 37--48.


\bibitem[Wahl1981]{Wahl 1981} J. Wahl, {\em Smoothings of normal surface singularities},
    Topology 20 (1981), 219--246.

 \bibitem[Xu-Yau1993]{Xu-Yau1993} Y.-J. Xu, S.S.-T. Yau, {\em Durfee conjecture and coordinate
     free characterization of homogeneous singularities.} J. Differential Geom. 37 (1993), no.
     2, 375--396.

\bibitem[Yau-Zhang2006]{Yau-Zhang2006} St.-T. Yau, L. Zhang, {\em An upper estimate of integral
    points in real simplices with an application to singularity theory.} Math. Res. Lett. 13
    (2006), no. 5--6, 911--921.
\end{thebibliography}
\end{document}